\documentclass[nohyperref]{article}

\usepackage{preamble}

\usepackage[accepted]{icml2022}

\usepackage[textsize=tiny]{todonotes}

\icmltitlerunning{Generalized Results for the Existence and Consistency of the MLE in the BTL Model}

\begin{document}

\twocolumn[
\icmltitle{Generalized Results for the Existence and Consistency of the MLE in the Bradley-Terry-Luce Model}



\icmlsetsymbol{equal}{*}

\begin{icmlauthorlist}
\icmlauthor{Heejong Bong}{cmu}
\icmlauthor{Alessandro Rinaldo}{cmu}
\end{icmlauthorlist}

\icmlaffiliation{cmu}{Department of Statiscis and Data Sciences, Carnegie Mellon University, Pittsburgh, PA, USA}

\icmlcorrespondingauthor{Heejong Bong}{hbong@andrew.cmu.edu}
\icmlcorrespondingauthor{Alessandro Rinaldo}{arinaldo@cmu.edu}

\icmlkeywords{Pairwise Comparison, Ranking, Bradley-Terry-Luce Model, Maximum Likelihood Estimator, Existence, Consistency, Machine Learning, ICML}

\vskip 0.3in
]



\printAffiliationsAndNotice{}  

\begin{abstract}
Ranking tasks based on pairwise comparisons, such as those arising in online gaming, often involve a large pool of items to order.  In these situations, the gap in performance between any two items can be significant, and the smallest and largest winning probabilities can be very close to zero or one. Furthermore, each item may be compared only to a subset of all the items, so that not all pairwise comparisons are observed. In this paper, we study the performance of the Bradley-Terry-Luce model for ranking from pairwise comparison data under more realistic settings than those considered in the literature so far. In particular, we allow for near-degenerate winning probabilities and arbitrary comparison designs. We obtain novel results about the existence of the  maximum likelihood estimator (MLE) and the corresponding $\ell_2$ estimation error without the bounded winning probability assumption commonly used in the literature and for arbitrary comparison graph topologies. Central to our approach is the reliance on the Fisher information matrix to express the dependence on the graph topology and the impact of the values of the winning probabilities on the estimation risk and on the conditions for the existence of the MLE. Our bounds recover existing results as special cases but are more broadly applicable.
\end{abstract}

\section{Introduction} \label{sec:introduction}
The task of eliciting a ranking over a possibly large collection of items based on partially observed pairwise comparisons is of great relevance for an increasingly broader range of applications, including sports analytics, multi-player online games, collection/analysis of crowdsourced data, marketing, aggregation of online search engine response queries and social choice models, to name a few. 
In response to these growing needs for theories and methods for high-dimensional ranking problems, the statistical and machine learning literature has witnessed a recent flurry of activities and novel results concerning high-dimensional parametric ranking models. Among them, the renown Bradley-Terry-Luce (BTL) model \citep{bradley1952rank,luce2012individual} has experienced a surge in popularity and use, due to its high interpretability and expressive power. 

Though the statistical literature on the BLT model is extensive, dating back to the seminal work on existence of the MLE by \cite{ford1957solution}, a series of recent contributions have led to accurate upper and lower bounds on the estimation error of the BTL model parameters under various high-dimensional settings. See Related Work below.
These bounds, as well as most other theoretical results in the literature, are typically predicated on the key assumption that the pairwise winning probabilities remain bounded away from $0$ and $1$ by a fixed and even  known amount. However, in high-dimensional ranking modeling, it is desirable to allow for the winning probabilities to become degenerate (i.e. to approach $0$ or $1$)  as the number of items to be compared grows. This feature expresses characteristics found in real-life ranking tasks, such as ratings in massive multi-player online games. In those cases, it is  expected that the maximal performance gap -- also known as dynamic range -- among different players increases with the number of players, so that the best players will be vastly superior to those with poor skills. Another significant complication frequently arising in real-life ranking tasks is that,  due to experimental  or financial constraints, only few among all possible pairwise comparisons are observed, resulting in a possibly irregular or sparse comparison graph. The combined effect of a broad gap in performance and of a possibly unfavorable comparison graph topology may have a considerable impact on the estimability of the BTL model parameters. The main goal of this paper  is to quantify such impact in novel and explicit ways.  We will expand on some of the latest theoretical contributions to parametric ranking using the BTL model in order to (i) formulate novel conditions for the existence of the MLE of the model parameters and (ii) derive more general and sharper $\ell_2$ estimation bounds. In particular, we will relax the prevailing assumption of non-degenerate winning probabilities frequently used in the literature on ranking while allowing for arbitrary comparison graph topologies. Our bounds hold under much broader settings than those considered so far and thus offers theoretical validation for the use of the BTL model in  real-life problems.

This paper makes the following two main contributions.

\begin{itemize}
    \item In \cref{sec:existence}, we present a novel sufficient condition for the existence of the MLE of the parameters of the BTL model under an arbitrary comparison graph and dynamic ranges that encompass and greatly generalize previous results. Existence of the MLE is a necessary pre-requisite for estimability of the model parameters, and therefore it is important to identify these minimal conditions for which the ranking task is at least well-posed. Compared to the few existing results in the literature, which hold only under very specific settings, our condition accounts for the effect of the graph topology and of the performance gap and, despite its generality, immediately recovers the few known results   as special cases; see \citet{simons1999asymptotics},\citet{yan2012sparse}, and \citet{han2020asymptotic}. 
    \item In \cref{sec:consistency}, we present a new $\ell_2$ consistency rate for the MLE without bounded dynamic range and again permitting arbitrary graph topologies that generalizes the $\ell_2$ error bound of \cite{shah2016estimation} and is more generally applicable. Two key technical contributions play a central role in our analysis:  (i) the use of a proxy loss function  that, unlike the negative log-likelihood loss, will retain strong convexity behavior even if the dynamic range is not bounded and (ii) the reliance of the Fisher information to adequately capture the efficiency of the MLE and its dependence on both the graph topology and the performance gap. In \Cref{sec:vs_chen_shah} we illustrated in detail how our new bound improves  the current state-of-the-art result of \cite{shah2016estimation} and is tighter in more general and challenging ranking settings. These findings are further corroborated in simulation studies. 
\end{itemize}

\paragraph{Recent Related Work.}
Over the last few years the high-dimensional statistics literature has produced an impressive body of new results on the  BTL model, demonstrating its performance and effectiveness in a variety of complex ranking tasks scenarios.
We present a non-exhaustive list of the related works. Novel bounds on the $\ell_2$ estimation of the BLT models are derived in \citet{negahban2012iterative,hajek2014minimax,khetan2016computational,shah2016estimation,chen2020partial}. \citet{hendrickx2019graph,hendrickx2020minimax} consider instead the sine error measure, while \citet{chen2019spectralregmletopk,chen2020partial,chen2021full} focus on the $\ell_\infty$ estimation error. 
To date, the sharpest $\ell_2$ related estimation rates are those obtained by \citet{shah2016estimation,hendrickx2020minimax}, and exhibit an explicit dependence on the topology of the comparison graph, i.e. the undirected graph over the items to be ranked, with edges indicating the pairs of items that are compared.

\paragraph{Notation.}
For symmetric matrices $A, B \in \reals^{d \times d}$, we write $A \succ B$ ($A \succeq B$) if $A - B$ is positive (semi-)definite. Also, for a vector $x \in \reals^d$, $\norm{x}_A^2 := x^\top A x$. We denote the eigenvalues of a positive semi-definite matrix $A$ by $0 \leq \lambda_1(A) \leq \lambda_2(A) \leq \dots \leq \lambda_d(A)$. We will use repeatedly the well known fact that if $A$ is the graph Laplacian or the Fisher information matrix induced by the BTL model then $A \mathbf{1} = \mathbf{0}$ and hence $\lambda_1(A) = 0$. Finally, for positive real sequences $x_d$ and $y_d$ with respect to $d$, we write
$
x_d = 
    o(y_d) ~\text{if}~ \frac{x_d}{y_d} \rightarrow 0,~
    O(y_d) ~\text{if}~ \frac{x_d}{y_d} ~\text{is upper bounded},~
    \omega(y_d) ~\text{if}~ \frac{x_d}{y_d} \rightarrow \infty, \textand
    \Omega(y_d) ~\text{if}~ \frac{x_d}{y_d} ~\text{is lower bounded}
$
as $d \rightarrow \infty$. In particular, we say $x_d = \Theta(y_d)$ if $x_d = O(y_d)$ and $x_d = \Omega(y_d)$. Accordingly, we write $x_d = o_p(y_d), O_p(y_d), \omega_p(y_d)$, or $\Omega_p(y_d)$ if the respective condition meets in probability; e.g., $x_d = O_p(y_d)$ if $\Pr[\frac{x_d}{y_d} \leq c] \rightarrow 0$ as $d \rightarrow \infty$ for some constant $c > 0$. 

\section{Bradley-Terry-Luce Model for Pairwise Comparisons} \label{sec:setting}

In this section, we review the Bradley-Terry-Luce model for ranking with pairwise comparisons introduced by \citet{bradley1952rank} and \citet{luce2012individual}. In this model, the $d$ items to be ranked   are each assigned an unknown numeric {\it quality scores}, denoted by the vector $w^* = (w^*_1,\ldots,w^*_d) \in \reals^d$ so that the probability that item $i$ is preferred to item $j$ in a pairwise comparison is specified as
\begin{equation}
    \Pr[i ~~\text{defeats}~~ j] = \frac{e^{w^*_i}}{e^{w^*_i}+e^{w^*_j}}.
\end{equation}
To ensure identifiability of the model, we impose the additional constraint that $\sum_{i=1}^d w^*_i = 0$.
The difference in performance between the best and worst items, that is, $\max_{i,j} \abs{w_i^* - w_j^*}$, is referred to as the maximal performance gap or {\it dynamic range}.

Data are collected in the form of $n$ pairwise comparisons among the items, which are modeled as a sequence of independent labeled Bernoulli random variables $(y_1,\ldots,y_n)$, where the label for $y_k$ is the pair of distinct items $(i_k, j_k) \in [d] \times [d]$ that have been compared. Accordingly, for  $k=1,\ldots,n$,
\begin{equation} \label{eq:BTL_model}
y_k = \left\{ 
\begin{array}{cc}
1  & \text{w.p. }\; \frac{\exp(w^*_{i_k})}{\exp(w^*_{i_k})+\exp(w^*_{j_k})}, \\
 -1  & \text{w.p. }\; \frac{\exp(w^*_{j_k})}{\exp(w^*_{i_k})+\exp(w^*_{j_k})}.
\end{array}
\right.
\end{equation}
Then, the random variable $X_k = y_k (\mathbf{e}_{i_k} - \mathbf{e}_{j_k}) \in \reals^d$ indicates the result of the $k$-th comparison between items $i_k$ and $j_k$: the item with entry $1$ won while the other with entry $-1$ was defeated. We record which items have been compared and the number of comparisons between any two items using the comparison matrix $X$, whose $k$-th row is $X_k^\top$. Then, the matrix
\begin{equation}
    L := \frac{1}{n} X^\top X.
\end{equation}
is the normalized Laplacian of the {\it comparison graph,} the weighted undirected graph over the items in which the weight of the $(i,j)$ edge  is the fraction of times, out of the  $n$ comparisons, that items $i$ and $j$ were compared. We note that $L$ turns out to be a function of $i_k$'s and $j_k$'s -- i.e., $L = \frac{1}{n} \sum_{k=1}^n (e_{i_k} - e_{j_k}) (e_{i_k} - e_{j_k})^\top$ -- so that $L$ is deterministic and also invariant with respect to the order between $i_k$ and $j_k$. 

We estimate the quality scores and then rank the items based on the maximum likelihood estimator. The (normalized) log-likelihood $l(w)$ of the BTL model is written in terms of the $X_k$'s as
\begin{equation} \label{eq:loglikelihood}
    l(w) = - \frac{1}{n} \sum_{k=1}^n \log(1+e^{-\inner{w,X_k}}),
\end{equation}
and the maximum likelihood estimator (MLE) $\hat{w}$ is obtained by solving the optimization problem
\begin{equation} \label{eq:mle}
    \hat{w} = {\arg\min}_{w \in \reals^d} -l(w) ~~~~~ \text{s.t.} ~ \sum_{i=1}^d w_i = 0.
\end{equation}
Notice that the the supremum in the above expression may not be attained, in which case we say the MLE does not exist. 

To obtain a tight consistency rate for the MLE, we need to take into account the both information about the topology of the comparison graph and the distribution of performance across items. As our results demonstrate, the Fisher information  matrix efficiently captures the impact of both sources of statistical difficulty. We note that under the BTL model, the Hessian matrix of the likelihood is constant with respect to $X$ and 
the Fisher information is easily calculated as 
\begin{equation}
\begin{split}
    & \mathcal{I}^*_{ij} := [\mathcal{I}(w^*)]_{ij} = \Exp[-\nabla^2 l(w^*)]_{ij} \\
    & = \begin{cases}
        \frac{L_{ij}}{\left(e^{(w^*_i - w^*_j)/2} + e^{(w^*_j - w^*_i)/2}\right)^2}, & i \neq j, \\
        -\sum_{k: k \neq i} [\mathcal{I}(w^*)]_{ik}, & i = j.
    \end{cases}
\end{split}
\label{eq:Fisher.info}
\end{equation}
One may easily recognize that the $(i,j)$-th entry of $\mathcal{I}(w^*)$ is equal to the corresponding entry of the normalized Laplacian $L$ weighted by $(e^{(w^*_i - w^*_j)/2} + e^{(w^*_j - w^*_i)/2})^{-2}$, which is decreasing in the performance gap $\abs{w^*_i - w^*_j}$. In the following sections, we will derive a novel; sufficient condition for the MLE existence and a bound on the magnitude of its error 
using the Fisher information $\mathcal{I}^*$ instead of the Laplacian $L$.


\section{Existence of the MLE under an Arbitrary Graph Topology} \label{sec:existence}

The MLE defined in \cref{eq:mle} often plays a fundamental role in any ranking task based on the BTL model. However, such an optimization-based estimator is not generally guaranteed to exist. 
For example, suppose that there is a partition of the items into two or more group  such that there is no comparison across different groups {\it by design.} Equivalently, the comparison graph is {\it disconnected}. In this case, no amount of data will be able to resolve the  ranking between any two items belonging to different groups, and the model is non-identifiable. As a result, the MLE does not exist. To avoid such cases, we will assume throughout that the comparison graph is fully connected, as is commonly done in the literature. 
Under this assumption, the Laplacian and Fisher information matrices have only one zero eigenvalue (with associated eigenvector spanning the one-dimensional subspace comprised of vectors with constant entries), and hence $\lambda_2(L)$ -- known as the {\it algebraic connectivity} of the graph \citep{brouwer2011spectra} --  and $\lambda_2(\mathcal{I}^*)$ are always positive.  

Existence and uniqueness of the solution of \cref{eq:mle} is still not guaranteed even with a connected comparison graph and demand stronger conditions. Indeed, when one item wins all comparisons against all the others, the log-likelihood function will not admit a supremum over $\mathbb{R}^n$ as the optimum will be achieved when a positive infinity score is assigned to the undefeated item. To avoid such cases, no single item nor group of items should overwhelm the others. As proved by \citet{Zermelo1929} and \citet{ford1957solution},  this is indeed sufficient and necessary for a  unique solution. We formalize this condition below. 

\begin{condition}[\citet{ford1957solution}] \label{cond:existence} 
In every possible partition of the items into two nonempty subsets, some item in the first set has beaten another in the second set. Or equivalently, for each ordered pair $(i, j)$, there exists a sequence of indices $i_0=i, i_1, \dots, i_m=j$ such that $i_{k-1}$ have beaten $i_k$ for $k=1,\dots,m$.
\end{condition}

\cref{cond:existence}  is informative about the existence of the MLE only after the results of all the comparisons are observed. An important yet not fully established problem is the prediction of \cref{cond:existence} ahead of observing any data. That is, we seek conditions for the score $w^*$ and the comparison graph, encoded by the Laplacian  $L$, to generate comparison results satisfying \cref{cond:existence} with high probability. To the best of our knowledge, this problem has been  studied in only few, limited settings and a general analytic result is lacking.
\citet{simons1999asymptotics} stated and proved a sufficient condition for the complete comparison graph, which demands a sufficient amount of comparisons and small dynamic range. This condition was further  extended by \cite{yan2012sparse} to allow for sparse graph topologies. More recently, \cite{han2020asymptotic} significantly sharpened Simon's and Yao's analysis and provided a sufficient condition for sparse Erd\"{o}s-Renyi type of graph topologies \citep{erdHos1960evolution}. 
In our first result, we further refine the arguments of \citet{simons1999asymptotics} to derive a novel and general sufficient condition for existence of the MLE that is valid for arbitrary comparison graphs and depends on the smallest non-zero eigenvalue of the Fisher information matrix. See \cref{app:pf_existence} for the proof.

\begin{theorem} \label{thm:existence}
If $\lambda_2(\mathcal{I}^*) \geq 2\frac{\log d}{n}$, then 
\begin{equation} \label{eq:existence}
    \Pr[\text{\cref{cond:existence} fails}] \leq \frac{2}{\sqrt{d}}.
\end{equation}
\end{theorem}

We emphasize that the above sufficient condition is not only generally applicable to any graph topology but in fact, as we will show next, recovers existing results as special cases. We believe that the strength and generality of our result stem from  the direct use of the Fisher information matrix.

\subsection{Comparison to Existing Results}

\noindent {\bf The complete graph.}
Under complete comparison graphs, \cref{thm:existence} recovers Lemma 1 of \citet{simons1999asymptotics}. In the setting, every pair of the items are compared the same number of times. Let $M := \max_{i,j \in [d]} \exp(w_i^* - w_j^*)$ --  a quantity analogue to the term $M_t$ in \citet{simons1999asymptotics} -- represent the dynamic range of the compared items. Then, according to \cref{eq:Fisher.info}, every off-diagonal element in $\mathcal{I}^*$ has absolute value larger than $\frac{M}{(1+M)^2}$, and it is easy to show that $\mathcal{I}^* \succeq \frac{M}{(1+M)^2}L$. Since the algebraic connectivity of the complete graph is given by $\lambda_2(L) = \frac{2}{d-1}$ \citep{brouwer2011spectra}, we arrive at the bound
\begin{equation} \label{eq:lambda_2_I_lambda_2_L}
    \lambda_2(\mathcal{I}^*) \geq \frac{M}{(1+M)^2} \lambda_2(L) \geq \frac{1}{1+M} \frac{1}{d-1}.
\end{equation}
As a result, the sufficient condition $M = o\left(\frac{d}{\log d}\right)$ of Lemma 1 of \cite{simons1999asymptotics} implies that $\lambda_2(\mathcal{I}^*) = \omega(\frac{\log d}{n} )$ which by \cref{thm:existence} implies existence of the MLE (to be precise, Lemma 1 of \cite{simons1999asymptotics} assumes the more restrictive scaling $M = o\left(\sqrt{\frac{d}{ \log d}}\right)$, but in fact their result holds under the weaker scaling reported above). 

\noindent {\bf The Erd\"os R\'enyi graph.}
As a second application, assume that the comparison graph is an  Erd\"os-R\'enyi graph with parameter $p$, possibly vanishing in $d$, a setting that has been the focus of much of the recent work on the BTL model \citep{negahban2012iterative,yan2012sparse,han2020asymptotic}. In this case, the sharpest known sufficient condition for the existence of the MLE of the BTL model parameters is due to \cite{han2020asymptotic}, and requires that $M = o(\frac{dp}{\log d})$, provided that $p = \omega(\frac{\log d}{d})$ and assuming again for simplicity a constant number of comparisons for each pair; see also \cite{yan2012sparse} for an earlier result requiring $p$ to be bounded away from $0$. Under the slightly weaker requirement that $p$ is of order  $\Omega (\frac{\log d}{d})$, the bound in  \eqref{eq:lambda_2_I_lambda_2_L} along with  the fact that
$\lambda_2(L) = \Theta_p(\frac{dp}{n})$ \citep[a fact that can be deduced from][after noting that we use a different normalization for the Laplacian]{kolokolnikov2014algebraic} shows that the same scaling of $M = o\left(\frac{dp}{\log d}\right)$ found by \citet{han2020asymptotic} does in fact satisfies the condition of \cref{thm:existence}.  


\section{Consistency of MLE without the Bounded Dynamic Range Assumption} \label{sec:consistency}

To the best of our knowledge, the work by \citet{shah2016estimation} was the first attempt to study  consistency of the MLE  of the BTL parameters in the $\ell_2$ and weighted $\ell_2$ loss  under arbitrary comparison graphs. Their rates are complemented by minimax lower bounds that also exhibit an explicit dependence on comparison graph topology, and can be used to identify comparison graph topologies  that produce  minimax optimal estimation rates. These results crucially require the assumption, frequently imposed in the literature on ranking,  of  a {\it bounded dynamic range}, namely that  $\norm{w^*}_\infty < B$, for a {\it known} value $B$. 
Under these constraint, the feasible space for the likelihood optimization problem \cref{eq:mle} is a compact set $\{ x \in \mathbb{R}^d \colon \| x\|_\infty \leq B\}$, over which the negative log-likelihood is by construction strongly convex. In particular, the resulting  regularized MLE is guaranteed to exist, even in cases in which the MLE may not. Next, using a Taylor series expansion, the authors obtain at a strongly quadratic approximation of the log-likelihood in a neighborhood of the true parameter value, from which the error upper bounds of their Theorems 1 and 2 follow.
However, the assumption imposed by \citet{shah2016estimation} is fairly restrictive and unrealistic in high-dimensional settings. 
Furthermore, the use of the constrained MLE is only feasible when prior information about the value of the hyper-parameter $B$ is available, a condition that is unrealistic in many settings. Hence, a consistency rate for the  MLE  that does not depend on the bounded dynamic range assumption is in order. 

Without the aid of the bounded dynamic range assumption, the summand $\log(1+e^{- \inner{w^*, X_k}})$ in the negative log-likelihood function $w \in \mathbb{R}^d \mapsto -l(w)$ gets increasingly flatter as the absolute value of the exponent increases, and hence the curvature of $-l(\cdot)$ is no longer uniformly lower-bounded.
In fact, the restricted curvature assumption is central to the analysis of \citet{shah2016estimation} and their techniques no longer apply when $-l(\cdot)$ is not strongly convex.
We overcome this hurdle by focusing instead on a proxy function $h_d : \reals^d \rightarrow \reals^d$ for the negative log-likelihood that offers two key advantages: (i) it turns $-l(\cdot)$ into a strongly convex loss function with respect to $h(\cdot - w^*)$ and (ii) it relates to the Fisher information matrix through the inequality
\begin{equation}
    \Exp[l(w^*) - l(w)] \geq h_d(w - w^*)^\top \mathcal{I}^* ~ h_d(w - w^*).
\end{equation}
With these critical modifications in place, the $\ell_2$ consistency rate of \cref{thm:consistency}, the main result of this section, follows using similar proof arguments as those deployed in \citet{shah2016estimation}. 

\begin{figure}[t!]
  \centering \setlength{\labelsep}{-3mm}
  \sidesubfloat[]{
    \includegraphics[height = 0.6\textwidth]{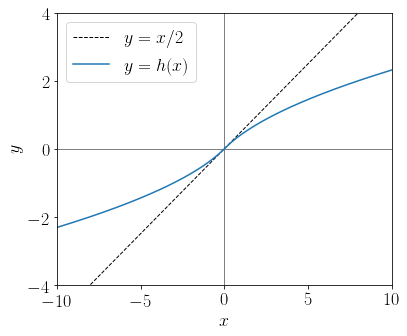}
    \label{fig:h_x}
  } \\
  \sidesubfloat[]{
    \includegraphics[height = 0.6\textwidth]{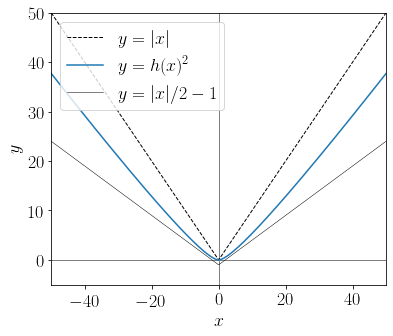}
    \label{fig:h2_x}
  }
  \caption{\sl {\bf Shape of the proxy function, $h$, and its square.} 
  {\bf (a)} 
  The absolute value of $h(x)$ is smaller than $\frac{\abs{x}}{2}$, and hence $h(x)^2$ is upper bounded by $\frac{x^2}{4}$. {\bf (b)} Also, $h(x)^2$ is upper- and lower bounded by $\abs{x}$ and $\frac{\abs{x}}{2 }- 1$, respectively. In sum, $h(x)^2$ behaves like a quadratic function around zero and an absolute function for sufficiently large $\abs{x}$.}
  \label{fig:h_and_h2}
\end{figure}

Before stating our main result, we first define and study the proxy function $h_d$. First, let $h:\reals \rightarrow \reals$ be the univariate real function defined by
$
  h(x) := \sgn(x) (\sqrt{\abs{x}+1}-1).
$
See \cref{fig:h_x}. We point out that $h$ is an increasing function that behaves linearly near the origin and whose slope decreases at the same rate as the slope of $\sqrt{|x|}$ for  $x \rightarrow \pm \infty$. Furthermore, it can be verified that the function $x \mapsto h(x)^2$ is upper-bounded and lower-bounded by $\abs{x}$ and $\frac{\abs{x}}{2}-1$, respectively; see  \cref{fig:h2_x}. As a result, $h(x)^2$ behaves like a quadratic function of $x$ near zero while  increasing linearly for sufficiently large $\abs{x}$. It is precisely this property that ensures strong convexity of $l(w)$ with respect to the proxy function $h_d(w - w^*)$, which we define next.
We take $h_d$ to be the $d$-dimensional element-wise evaluation of $h$. That is, 
$
    [h_d(x)]_i = h(x_i), ~~ i = 1, \dots, d.
$
As remarked above, when  $x \rightarrow \infty$, $\norm{h_d(x)}_2^2 = \sum_{i=1}^d h(x_i)^2$ behaves like a linear function, so its curvature converges to zero. Hence, although $w \in \mathbb{R}^d \mapsto -l(w)$ is not itself strongly convex, the MLE problem in \cref{eq:mle} can be approximated by a quadratic problem with respect to $w \in \mathbb{R}^d \mapsto h_d(w - w^*)$. 

The $\ell_2$ convergence rate of the MLE derived below depends on a parameter $\kappa$ 
defined by
\begin{equation}
\begin{split}
    \kappa := & \lambda_{d} (L^{\frac{1}{2}} \mathcal{I}^{*+} L^{\frac{1}{2}}) 
\end{split}
\end{equation} 
where $\mathcal{I}^{*+}$ is the Moore-Penrose pseudo-inverse of the Fisher information matrix. As a function of $\mathcal{I}^*$, $\kappa$ takes into account  both the comparison graph topology and the gaps in performance among the items. Importantly, $\kappa$ is determined  only by  the performance gaps among compared items, unlike the the dynamic range, which is the maximal performance gap across {\it all} items, including those that are never compared. We will elucidate the gains from replacing the traditional dynamic range parameter by the new parameter $\kappa$ through a concrete example below, when we compare our bounds to existing results.  

\begin{theorem} \label{thm:consistency} 
Assume that $\lambda_2(\mathcal{I}^*) \geq c_0 \frac{\sqrt{v_\text{max} \log d}}{n}$ for some universal constant $c_0 > 0$, where $v_\text{max}$ be the maximum degree of the comparison graph. Then, the MLE $\hat{w}$ for the BTL model exists with probability at least $1 - \frac{1}{d}$, and the error $\Delta = \hat{w} - w^*$ satisfies
\begin{equation}\label{eq:ell_I}
    \norm{h_d(\Delta)}_{\mathcal{I}^*}^2 \leq c_1 \kappa \frac{td}{n},
\end{equation}
and, as a corollary,
\begin{equation}\label{eq:ell_2_alt}
    \norm{h_d(\Delta)}_2^2 \leq c_1 \frac{\kappa}{\lambda_2(\mathcal{I}^*)} \frac{td}{n},
\end{equation}
with probability at least $1 - e^{-t} - \frac{1}{d}$ for any $t>0$ and some universal constant $C > 0$.
\end{theorem}


The only assumption in the theorem is the requirement that
$\lambda_2(\mathcal{I}^*) \geq c_0 \frac{\sqrt{v_\text{max} \log d}}{n}$, which is imposed to ensure existence of the MLE and a sufficient degree of curvature of the log-likelihood. As remarked above, this issue does not arise with the regularized MLE, which is always well defined. 

We point out that the condition $\lambda_2(\mathcal{I}^*) \geq c_0 \frac{\sqrt{ v_\text{max} \log d}}{n} $ is in general stronger than the sufficient condition for the existence of the MLE implied by \Cref{thm:existence}. To see this, we use the well known fact that the maximal degree is of the same order as the largest eigenvalue of the graph Laplacian, which in our settings implies that $ v_\text{max}  = \Theta(\lambda_2(\mathcal{I}^*) n)$. Thus,  $\lambda_2(\mathcal{I}^*) \geq c_0 \frac{\sqrt{  v_\text{max} \log d}}{n} $ is equivalent to $\lambda_2(\mathcal{I}^*) = \Omega( \frac{\lambda_{d}(\mathcal{I}^*)}{\lambda_2(\mathcal{I}^*)}\frac{ \log d}{n})$. Ignoring constants, the latter condition is stricter than the condition $\lambda_2(\mathcal{I}^*) \geq 2\frac{\log d}{n}$ from \Cref{thm:existence}, unless all the non-zero eigenvalues of $\mathcal{I}^*$
are of the same order. It is an open problem to fill the gap between the conditions for the MLE existence and $\ell_2$ error rate guarantee.

Below we provide detailed comparisons to existing results and explain how \Cref{thm:consistency} delivers significant improvements by allowing for a significant weakening of the assumptions.

\subsection{Comparison to Existing Results} \label{sec:vs_chen_shah}

\noindent {\bf Erd\"os-R\'enyi graph under bounded dynamic range assumption.}
We first compare \Cref{thm:consistency} to Theorem 3.1 of \citet{chen2019spectralregmletopk}, who provide an $\ell_2$ rate for the BTL parameters assuming an Erd\"os-R\'enyi comparison graph with edge probability $p$ and bounded dynamic range $B$. Theorem 3.1 of \citet{chen2019spectralregmletopk} requires that $p = \Omega(\frac{\log d}{d})$, to ensure that the comparison graph is connected. This additional constraint is equivalent to the condition $\lambda_2(\mathcal{I}^*) \geq c_0 \frac{\sqrt{ v_{\max} \log d}}{n}$  of \Cref{thm:consistency}. According to the well-known probabilistic guarantees for the Laplacian of the adjacency matrix of an ER model, $\lambda_2(\mathcal{I}^*) = \Theta_p(\frac{dp}{n})$ and $v_\mathrm{max}=\Theta_p(dp)$ \citep[see, e.g.,][]{MAL-048}. Thus, with high probability, 
\begin{equation}
\begin{split}
    \lambda_2(\mathcal{I}^*) \geq c_0 \frac{\sqrt{v_{\max} \log d}}{n} 
    \Leftrightarrow & \frac{dp}{n} = \Omega\left(\frac{\sqrt{dp \log d}}{n}\right) \\
    \Leftrightarrow & p = \Omega\left(\frac{\log d}{d}\right),
\end{split}
\end{equation}
for some universal constant $c_0 > 0$.
We mention in passing that, in this case, the requirements on $\lambda_2(\mathcal{I}^*)$ posed in \Cref{thm:consistency} and \Cref{thm:existence} are equivalent.



As for the $\ell_2$ error rates, our bound and the one in \citet{chen2019spectralregmletopk} are as follows:
\begin{equation} \label{eq:vs_chen}
\begin{array}{cc}
    \text{\Cref{thm:consistency}:} &  \norm{h_d(\Delta)}_2^2 \leq c_1 \frac{\kappa^2t}{pL}\\
    \text{\citet{chen2019spectralregmletopk}:} & \norm{\Delta}_2^2 \leq c' \frac{1}{pL}
\end{array}
\end{equation}
with probability at least $1-e^{-t}-\frac{1}{d}$, for any $t > 0$, and $1-O(d^{-7})$, respectively. Above, $c'$ and $\kappa$ are quantities depending on the dynamic range $B$. In the Erd\"os-R\'enyi graph case, $\kappa \leq (e^{-B} + e^{B})^2$, where the dependence of $c'$ on the dynamic range is implicit. The results in \cref{eq:vs_chen} are similar except for the use of the proxy function in \Cref{thm:consistency}. When the right hand side of the first inequality in \cref{eq:vs_chen} is bounded, $h_d(\Delta)$ can be lower bounded by $\Omega(\|\Delta\|_2^2) $, so that the both theorems yield essentially the same $\ell_2$ convergence rate. 
The fact that our general bound, applicable to arbitrary comparison graphs, yields the same convergence rate as Theorem 3.1 in \citet{chen2019spectralregmletopk} is remarkable. Indeed, the analysis of \citet{chen2019spectralregmletopk} is specifically tailored to the Erd\"os-R\'enyi graph, as the authors relied crucially on the special properties and the high degree of regularity of  this type of graph. 


\noindent {\bf Arbitrary graph  under bounded dynamic range assumption.}
%
%
%
We now show that the error bound of \Cref{thm:consistency} is similar to the one derived in \citet{shah2016estimation} under the same settings considered by the authors, namely an arbitrary comparison graph and bounded dynamic range $B$. It is worth repeating that our results concern the MLE of the BTL parameters and does not require knowledge of $B$, while \citet{shah2016estimation} make the additional convenient -- and arguably impractical -- assumption of a known, bounded dynamic range, which allows them to focus on the simpler constrained MLE. 

In detail, Theorem 1 of \citet{shah2016estimation} analyze the $\ell_\infty$ regularized MLE
\begin{equation}
    \tilde{w} := \arg\min_{w:\norm{w}_\infty \leq B} l(w)
\end{equation}
and provides the $\ell_2$-norm bound of the error $\tilde{\Delta} := \tilde{w}-w^*$ as
\begin{equation} \label{eq:ell_2_shah}
    \norm{\tilde{\Delta}}_{L}^2 \leq c'' \frac{\zeta^2}{\gamma^2}\frac{td}{n},
\end{equation}
with probability at least $1-e^{-t}$ for any $t > 0 $ where $c''$ is an universal constant. Here, $\gamma$ and $\zeta$ are parameters accounting for the curvature of the log-likelihood function and performance distribution across the items, defined by
\begin{equation}\label{eq:gamma.zeta}
\begin{array}{c}
    \gamma :=  \min_{x \in [-2B,2B]} \frac{d^2}{dt^2}(-\log\psi(x)) \\
    \zeta := \max_{x \in [0,2B]} \frac{\psi'(x)}{\psi(2x)(1-\psi(2x))}, \\
\end{array}
\end{equation}
and relate to the parameter $\kappa$ in \cref{eq:ell_2_alt}. Under the Bradley-Terry-Luce model where $\psi(x) = \frac{1}{1+e^{-x}}$, they are both functions of the dynamic range parameter $B$: $\gamma = \frac{1}{(e^{-B}+e^B)^2}$ and $\zeta = 1$.
Then, \cref{eq:ell_2_shah} becomes
\begin{equation} \label{eq:ell_2_bounded_shah}
    \norm{\tilde{\Delta}}_L^2 \leq c'' {(e^{-B}+e^B)^4}\frac{td}{n}.
\end{equation}

On the other hand, repeating the same argument leading up to \cref{eq:lambda_2_I_lambda_2_L}, we obtain that $\mathcal{I}^* \succeq \frac{M}{(1+M)^2} L$ where $M = e^{2B}$. Thus, under the bounded dynamic range, 
$
    \norm{h_d(\Delta)}_{\mathcal{I}^*}^2 \geq \frac{1}{(e^{-B}+e^B)^2} \norm{h_d(\Delta)}_{L}^2
$
and also $\kappa \leq (e^{-B}+e^B)^2$. 
Thus,  \cref{thm:consistency} yields that 
\begin{equation} \label{eq:ell_2_bounded_ours}
    \norm{h_d(\Delta)}_L^2 \leq c_1 {(e^{-B}+e^B)^4}\frac{td}{n}
\end{equation}
with probability at least $1 - e^{-t} - \frac{1}{d}$ for any $t > 0$. The bounds \eqref{eq:ell_2_bounded_shah} and  \eqref{eq:ell_2_bounded_ours} for the MLE and regularized MLE respectively (the latter requiring knowledge of the parameter $B$)  are identical except for a potential difference in the universal constants and the use of the proxy function in \cref{eq:ell_2_bounded_ours}. As explained in the previous paragraph, the proxy function yields the same convergence rate as in \cref{eq:ell_2_bounded_shah} if the left hand side of \cref{eq:ell_2_bounded_ours} is bounded; 
otherwise, \Cref{thm:consistency} provides a looser bound than \citet{shah2016estimation}. 

We also note that \citet{hendrickx2019graph} established a BTL parameter estimator with minimax property in the sine error measure. Under bounded dynamic range, the sine error measure is equivalent to $\ell_2$ error.
Hence, the $\ell_2$ estimation rate implied by \citet{hendrickx2019graph} is strictly sharper than one implied by ours and \citet{shah2016estimation}, which is expected. For example, in Section 4.1 therein, \citet{shah2016estimation} discussed the suboptimality of the $\ell_2$ upperbound in their Theorem 2 under particular types of comparison graphs because the upperbound does not match its companion lowerbound. Instead, the sharper rate of the lowerbound is achieved by the $\ell_2$ rate implication of Theorem 1, \citet{hendrickx2019graph}; compare, e.g.,  the examples in Section 4.1, \citet{shah2016estimation} and Table 1, \citet{hendrickx2019graph}. For the proper comparison, the sine error rate in the table should be squared and multiplied by ``$\frac{n^2}{k}$'' in their notation. (Note that ``$n$'' is the number of the compared items and ``$k$'' is the number of comparison in \citet{hendrickx2019graph}.) On the other hand, the estimation rates of \citet{shah2016estimation} are minimax in the $\norm{\cdot}_L$ error. Because our bound inherits the minimax property of \citet{shah2016estimation} under bounded dynamic range, it is generally better than \citet{hendrickx2019graph} in the $\norm{\cdot}_L$ error.


In the present discussion, since we are assuming a general comparison graph, we have bounded both the ratio $\zeta/\gamma$ and $\kappa$, appearing in \eqref{eq:ell_2_bounded_shah} and  \eqref{eq:ell_2_bounded_ours} respectively, using $B$. While $\zeta/\gamma$  are always functions of $B$ only regardless of the comparison graph, the parameter $\kappa$ also depends on the specific comparison graph topology. This dependence can in turn be leveraged in some specific cases to reduce the value of $\kappa$ and thus produce shaper rates via \cref{thm:consistency}. We elaborate on this point next.

\noindent{\bf Banded graph  without bounded dynamic range assumption.}
%
Theorem 4.2 of \citet{shah2016estimation} and \cref{thm:consistency} differ crucially in the ways they account in the error rate for the dependence on the gap in performance among items and on the comparison graph topology.
Theorem 4.2 of \citet{shah2016estimation} decouples these two aspects. The gap in performance enters the $\ell_2$ bound \eqref{eq:ell_2_bounded_shah} only through the ratio $\zeta/\gamma$, which depends on the dynamic range parameter $B$ but not the comparison graph; on the other hand, the graph topology affects the rate only through the algebraic connectivity $\lambda_2(L)$, which is independent of the gap in performance. In contrast, the $\ell_2$ rate in \cref{thm:consistency} depends on $\kappa$ and $\lambda_2(\mathcal{I}^*)$, both of which simultaneously quantify the impact of the graph topology and the gap in performance, and their potential interaction. When we crudely decouple such dependence, as we did above, such interaction is lost and we essentially recover the bound in Theorem 4.2 of \citet{shah2016estimation}. But a more careful use of the dependence will lead to better results by \cref{thm:consistency}. We illustrate this fact using a banded graph topology.

Suppose that $w^*_i$'s are evenly distributed along the dynamic range
\begin{equation}\label{eq:even_dist_players}
    w^*_i := \frac{2i-d}{d}B, ~~ i \in [d]
\end{equation}
and that 
the comparisons are made $T$ times only between items having a difference in ranking smaller than a fixed positive numbers $W > 0$. The resulting normalized graph Laplacian is 
\begin{equation} \label{eq:laplacian}
    L_{ij} = \begin{cases}
    -\sum_{k: k \neq i} L_{ik}, & i = j, \\
    -1/n, & 0 < \abs{i - j} \leq W, \\
    0, & \text{elsewhere,}
    \end{cases}
\end{equation}
and the total number of comparisons $n$ is $T(dW-W(W-1)/2)$. Since the graph Laplacian is a banded matrix, we call the graph a {\it banded comparison graph} and the hyperparameter $W$ the {\it comparison width}. Banded comparison graphs are easily seen in real data, in particular, when the size of the item pool is huge. In the high-dimensional setting, it is impossible to compare every pair of the items due to financial and physical constraints. To maximize the efficiency of limited resource, comparisons are often made only between the items of which none are overwhelming the other. The influence of this comparison graph to the estimation error is measured by its algebraic connectivity, which we establish in the following lemma. See \cref{app:pf_lambda2L} for the proof.
\begin{lemma} \label{lem:lambda2L}
    The algebraic connectivity of the given banded comparison graph is
    $
        \lambda_2(L) 
        = \Theta\left(\frac{W^2}{d^3}\right)
    $.
\end{lemma}

Since $\gamma$ and $\zeta$ for Theorem 2 in \citet{shah2016estimation} depend only on the dynamic range, their values are still as in \eqref{eq:gamma.zeta}, and the resulting $\ell_2$ error bound from Theorem 4.2 of \citet{shah2016estimation} is
\begin{equation}
\begin{split}
    \norm{\tilde{\Delta}}_2^2 \leq c'' \frac{\zeta^2}{\gamma^2\lambda_2(L)} \frac{td}{n} 
    & \leq 2 c'' \frac{td^3}{TW^3} e^{4B},
\end{split}
\end{equation}
with  probability at least $1 - e^{-t}$ for any $t > 0$ and for some the universal constant $c''$ (See \cref{eq:ell_2_bounded_shah}).

On the other hand, \cref{thm:consistency} leads to a tighter bound. To see this, the Fisher information matrix for this graph topology is
\begin{equation}
    \mathcal{I}^*_{ij} = \begin{cases}
    -\sum_{k: k \neq i} \mathcal{I}^*_{ij}, & i = j, \\
    -\frac{1/(d-1)}{(e^{-\frac{B}{d}\abs{i-j}}+e^{\frac{B}{d}\abs{i-j}})^2}, & 0 < \abs{i - j} \leq W, \\
    0, & \text{elsewhere,}
    \end{cases}
\end{equation}
where $\abs{\mathcal{I}_{ij}(w^*)} \geq \frac{1}{4e^{2BW/d}} \abs{L_{ij}}$ for every $(i, j)$ such that $0 < \abs{i - j} \leq W$. As a result, $\kappa$ and $\lambda_2(\mathcal{I}^*)$ turns out to be much smaller than in the previous example: $\kappa \leq 4e^{2BW/d}$ and $\lambda_2(\mathcal{I}^*) = \Omega\left(\frac{1}{d^3e^{2BW/d}}\right)$.
Putting all these pieces together, we arrive at a sharper $\ell_2$ error bound 
\begin{equation} \label{eq:ell_2_banded_ours}
\begin{split}
    \norm{h_d(\Delta)}_2^2 
    \leq & c_1 \frac{\kappa}{\lambda_2(\mathcal{I}^*)} \frac{td}{n}
    \leq c_1 \frac{td^3}{TW^3} e^{\frac{4BW}{d}}
\end{split}
\end{equation}
with probability at least $1 - \frac{1}{d} - e^{-t}$. 

When $B$ increases in $d$, \cref{eq:ell_2_banded_ours} provides a much tighter convergence rate so the $\ell_2$ convergence requires a much smaller sample complexity in terms of $d$. Concretely, suppose that $\norm{w^*}_\infty$ increases at a rate $\sqrt{\log d}$ -- this would be for example the case if the $w_i^*$'s were independently  sampled from a standard Gaussian distribution, with high probability. 
Then, Theorem 2 of \cite{shah2016estimation} implies that every compared pair should be matched $T = \omega(e^{8\sqrt{\log d}} \log d)$ times for $\ell_2$ consistency even if $W = d$. On the other hand, \cref{thm:consistency} suggests a much milder condition that every compared pair should be matched $T = \omega((\log d)^{5/2})$ times when $W = \Theta(\frac{d}{\sqrt{\log d}})$.

We remark that this result does not contradict the minimax lower bound found in Theorem 1 of \citet{shah2016estimation}. In fact, although their upper bound is proved to be optimal with respect to $d$ and $n$ under fixed $B$, it has a sub-optimal dependence on $B$ at least when compared with the minimax lower bound therein. 
In the above example,  \cref{thm:consistency} yields improvements in the dependence to $B$ and graph topology by incorporating them into Fisher information matrix $\mathcal{I}^*$.
In \cref{sec:simulation}, we demonstrate the improved dependence in \cref{thm:consistency} by simulations, particularly with respect to the dynamic range and graph topology.

\section{Illustrative Simulations} \label{sec:simulation}

\begin{figure}[t!]
  \centering 
  \setlength{\labelsep}{-1mm}
  \sidesubfloat[]{
    \includegraphics[height = 0.45\textwidth]{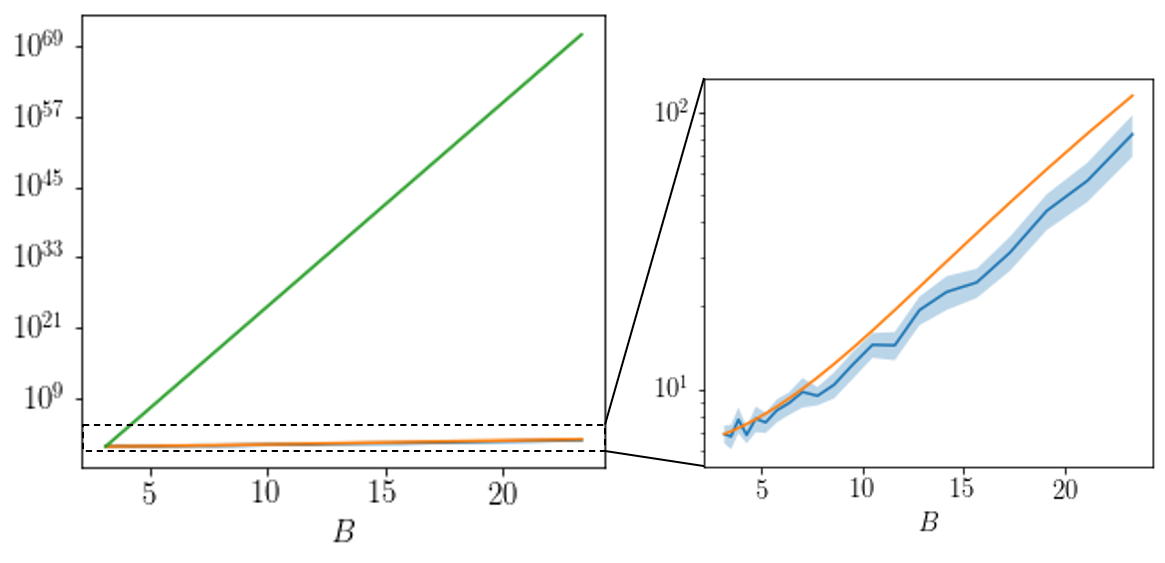}
    \label{fig:l2_vs_B}
  }
  \setlength{\labelsep}{-4mm} 
  \sidesubfloat[]{
    \includegraphics[height = 0.45\textwidth]{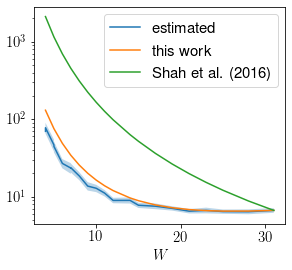}
    \label{fig:l2_vs_W}
  }

  \caption{\sl {\bf Estimated $\ell_2$ error by simulation compared to the theorems.} $\ell_2$ rates are estimated under various {\bf (a)} maximal performance differences $B$ and {\bf (b)} graph topologies induced by comparison widths $W$. With log-scaled y-axes, the blue lines indicate the averages of the $\ell_2$ errors over the $100$ simulation runs where the blue shaded areas represent the 95\% pointwise confidence interval for the averages. The yellow and green lines indicate the upperbounds of the $\ell_2$ errors suggested by \cref{thm:consistency} and Theorem 2 of \citet{shah2016estimation}, respectively, up to constant scales (which are constant shifts in the plots). The zoom pane only shows the estimated values and the suggested upperbound by \cref{thm:consistency} for the readability. With respect to both $B$ and $W$, \cref{thm:consistency} gives tighter bounds to the $\ell_2$ error.
  }
  \label{fig:l2_results}
\end{figure}

\citet{shah2016estimation} demonstrated the tightness of the error bounds in their Theorem 2 therein with respect to the number of the comparisons $n$ through simulations. We do the same for \cref{thm:consistency} but with respect to the dynamic range and the graph topology. The comparison data were simulated using  the BTL model (\cref{eq:BTL_model}) under the banded comparison graphs introduced in \cref{sec:vs_chen_shah}. In our experiment, the number of compared items $d$ and the number of comparisons between compared pairs $T$ are fixed to $100$ and $5$, respectively. The maximal performance gap $B$  ranges from $4 e^{-1} \sqrt{\log d}$ to $4 e \sqrt{\log d}$, and the comparison width $W$  from $e^{-1} \frac{d}{4 \sqrt{\log d}}$ to $e \frac{d}{4 \sqrt{\log d}}$. The MLE was calculated by iterative Luce Spectral Ranking \citep{maystre2015fast}, provided by Python package \texttt{choix}\footnote{Published through PyPi under MIT License}. For each pair of values $(B, W)$, we performed $100$ simulations, and the average of the estimated $\ell_2$ error is plotted in \cref{fig:l2_results} (blue line) along with the consistency rates implied by \cref{thm:consistency} (yellow line) and Theorem 2 of \citet{shah2016estimation} (green line). The light blue shaded area represents 95\% pointwise confidence intervals for the estimated averages. \cref{fig:l2_vs_B} shows the simulation results when $W$ was fixed at $\frac{d}{\sqrt{\log d}}$ and $B$ changed across the range while \cref{fig:l2_vs_W} displays the result for $B$ at the fixed value of $\sqrt{\log d}$ and $W$ changing. Because \cref{thm:consistency} and Theorem 2 of \citet{shah2016estimation} both provide $\ell_2$ consistency rates up to a universal scale, the yellow and green lines in \cref{fig:l2_results} are arbitrarily shifted to make the comparison easier. Hence, when comparing our bound (orange line) with that of \cite{shah2016estimation} (green line), we should not compare the relative position of the corresponding curves but their rates of increase as $B$ increases or $W$ decreases, especially in relation to the blue curve, based on the data. The plots clearly indicate that \cref{thm:consistency} provides tighter rates in both $B$ and $W$. We obtained similar results for varying  $B$ in other simulations involving different graph topologies. See \cref{app:supp_simulation} and \url{github.com/HeejongBong/mmpc} for the supplementary simulation results and reproducible simulation code scripts.

\section{Discussion} \label{sec:discussion}

We have derived novel conditions for the existence of the MLE and novel $\ell_2$ estimation rates in the BTL model under weaker assumptions than those assumed in the literature so far. In particular, we allow arbitrary comparison graphs and do not require a bounded gap among the quality scores, nor knowledge of such bound. 
Our bounds, based on the Fisher information matrix, 
 recover existing results as corollaries
and are also applicable to more general settings of practical relevance. 
Our theoretical results support the use of the BTL model for ranking problems and more generally statistical inference with pairwise comparison data. 

To the best of our knowledge, the use of proxy functions in analyzing the MLE solutions and their properties is new, although convex surrogate loss functions are widely used in classification. Because the BTL model is a special case of the logistic regression, one might use the same technique to establish a tight uncertainty measure for general logistic models when the odds are not necessarily bounded. For example, a promising open problem is the sine error minimax estimate of the BTL model without bounded dynamic range by incorporating the proxy function technique and \citet{hendrickx2019graph}.

Another important topic in the ranking literature which we did not  discuss here is the probability bound of identifying the $K$ best items from the observed partial orderings: namely, top-$K$ ranking. Recently, the minimax bound under the Erd\"os-R\'enyi comparison graph was established by \citet{chen2020partial,chen2019spectralregmletopk}, and a decision-theroetic recovery of full ranking was suggested by \citet{chen2021full}. However, the learning theory of top-$K$ ranking from passive observation of pair-wise comparisons in arbitrary graphs is yet unknown. As discussed in \citet{chen2020partial}, theoretical guarantee for top-$K$ recovery requires the $\ell_\infty$ estimation error to be controlled. To establish such bound, the performance score $w^*$ and comparison graph $L$ should be more fine-tuned. 


One might point out that 
the banded comparison graphs, used for the comparison against \citet{shah2016estimation}, already contained prior knowledge about $w^*$, which is not properly evaluated in \cref{sec:vs_chen_shah}. However, such comparison graphs can emerge without prior knowledge in online settings where the comparison graph is decided based on the recent estimate of the BTL parameters. A promising direction of future study is to study the theoretical property of such online estimators.



\section*{Acknowledgements}

We gratefully acknowledge helpful comments and mindful suggestions from the reviewers and meta-reviewer. We also send the deepest appreciation for inspiring discussions with Wanshan Li and Shamindra Shrotriya in the Department of Statistics and Data Sciences, Carnegie Mellon University. 




\bibliographystyle{icml2022}
\bibliography{references}  

\begin{thebibliography}{23}
\providecommand{\natexlab}[1]{#1}
\providecommand{\url}[1]{\texttt{#1}}
\expandafter\ifx\csname urlstyle\endcsname\relax
  \providecommand{\doi}[1]{doi: #1}\else
  \providecommand{\doi}{doi: \begingroup \urlstyle{rm}\Url}\fi

\bibitem[Bradley \& Terry(1952)Bradley and Terry]{bradley1952rank}
Bradley, R.~A. and Terry, M.~E.
\newblock Rank analysis of incomplete block designs: I. the method of paired
  comparisons.
\newblock \emph{Biometrika}, 39\penalty0 (3/4):\penalty0 324--345, 1952.

\bibitem[Brouwer \& Haemers(2011)Brouwer and Haemers]{brouwer2011spectra}
Brouwer, A.~E. and Haemers, W.~H.
\newblock \emph{Spectra of graphs}.
\newblock Springer Science \& Business Media, 2011.

\bibitem[Chen et~al.(2020{\natexlab{a}})Chen, Gao, and Zhang]{chen2020partial}
Chen, P., Gao, C., and Zhang, A.~Y.
\newblock Partial recovery for top-$ k $ ranking: Optimality of mle and
  sub-optimality of spectral method.
\newblock \emph{arXiv preprint arXiv:2006.16485}, 2020{\natexlab{a}}.

\bibitem[Chen et~al.(2020{\natexlab{b}})Chen, Gao, and Zhang]{chen2021full}
Chen, P., Gao, C., and Zhang, A.~Y.
\newblock Optimal full ranking from pairwise comparisons.
\newblock \emph{arXiv preprint arXiv:2101.08421}, 2020{\natexlab{b}}.

\bibitem[Chen et~al.(2019)Chen, Fan, Ma, and Wang]{chen2019spectralregmletopk}
Chen, Y., Fan, J., Ma, C., and Wang, K.
\newblock Spectral method and regularized {MLE} are both optimal for top-{$K$}
  ranking.
\newblock \emph{Ann. Statist.}, 47\penalty0 (4):\penalty0 2204--2235, 2019.
\newblock ISSN 0090-5364.
\newblock \doi{10.1214/18-AOS1745}.

\bibitem[Erd{\H{o}}s \& R{\'e}nyi(1960)Erd{\H{o}}s and
  R{\'e}nyi]{erdHos1960evolution}
Erd{\H{o}}s, P. and R{\'e}nyi, A.
\newblock On the evolution of random graphs.
\newblock \emph{Publ. Math. Inst. Hung. Acad. Sci}, 5\penalty0 (1):\penalty0
  17--60, 1960.

\bibitem[Ford~Jr(1957)]{ford1957solution}
Ford~Jr, L.~R.
\newblock Solution of a ranking problem from binary comparisons.
\newblock \emph{The American Mathematical Monthly}, 64\penalty0 (8P2):\penalty0
  28--33, 1957.

\bibitem[Hajek et~al.(2014)Hajek, Oh, and Xu]{hajek2014minimax}
Hajek, B., Oh, S., and Xu, J.
\newblock Minimax-optimal inference from partial rankings.
\newblock \emph{Advances in Neural Information Processing Systems}, 27, 2014.

\bibitem[Han et~al.(2020)Han, Ye, Tan, Chen, et~al.]{han2020asymptotic}
Han, R., Ye, R., Tan, C., Chen, K., et~al.
\newblock Asymptotic theory of sparse bradley--terry model.
\newblock \emph{Annals of Applied Probability}, 30\penalty0 (5):\penalty0
  2491--2515, 2020.

\bibitem[Hanson \& Wright(1971)Hanson and Wright]{hanson1971bound}
Hanson, D.~L. and Wright, F.~T.
\newblock A bound on tail probabilities for quadratic forms in independent
  random variables.
\newblock \emph{The Annals of Mathematical Statistics}, 42\penalty0
  (3):\penalty0 1079--1083, 1971.

\bibitem[Hendrickx et~al.(2019)Hendrickx, Olshevsky, and
  Saligrama]{hendrickx2019graph}
Hendrickx, J., Olshevsky, A., and Saligrama, V.
\newblock Graph resistance and learning from pairwise comparisons.
\newblock In \emph{International Conference on Machine Learning}, pp.\
  2702--2711. PMLR, 2019.

\bibitem[Hendrickx et~al.(2020)Hendrickx, Olshevsky, and
  Saligrama]{hendrickx2020minimax}
Hendrickx, J., Olshevsky, A., and Saligrama, V.
\newblock Minimax rate for learning from pairwise comparisons in the btl model.
\newblock In \emph{International Conference on Machine Learning}, pp.\
  4193--4202. PMLR, 2020.

\bibitem[Khetan \& Oh(2016)Khetan and Oh]{khetan2016computational}
Khetan, A. and Oh, S.
\newblock Computational and statistical tradeoffs in learning to rank.
\newblock \emph{Advances in Neural Information Processing Systems}, 29, 2016.

\bibitem[Kolokolnikov et~al.(2014)Kolokolnikov, Osting, and
  Von~Brecht]{kolokolnikov2014algebraic}
Kolokolnikov, T., Osting, B., and Von~Brecht, J.
\newblock Algebraic connectivity of erd{\"o}s-r{\'e}nyi graphs near the
  connectivity threshold.
\newblock 2014.

\bibitem[Luce(2012)]{luce2012individual}
Luce, R.~D.
\newblock \emph{Individual choice behavior: A theoretical analysis}.
\newblock Courier Corporation, 2012.

\bibitem[Maystre \& Grossglauser(2015)Maystre and
  Grossglauser]{maystre2015fast}
Maystre, L. and Grossglauser, M.
\newblock Fast and accurate inference of plackett-luce models.
\newblock Technical report, 2015.

\bibitem[Negahban et~al.(2012)Negahban, Oh, and Shah]{negahban2012iterative}
Negahban, S., Oh, S., and Shah, D.
\newblock Iterative ranking from pair-wise comparisons.
\newblock \emph{Advances in neural information processing systems},
  25:\penalty0 2474--2482, 2012.

\bibitem[Shah et~al.(2016)Shah, Balakrishnan, Bradley, Parekh, Ramchandran, and
  Wainwright]{shah2016estimation}
Shah, N.~B., Balakrishnan, S., Bradley, J., Parekh, A., Ramchandran, K., and
  Wainwright, M.~J.
\newblock Estimation from pairwise comparisons: Sharp minimax bounds with
  topology dependence.
\newblock \emph{The Journal of Machine Learning Research}, 17\penalty0
  (1):\penalty0 2049--2095, 2016.

\bibitem[Simons \& Yao(1999)Simons and Yao]{simons1999asymptotics}
Simons, G. and Yao, Y.-C.
\newblock Asymptotics when the number of parameters tends to infinity in the
  bradley-terry model for paired comparisons.
\newblock \emph{The Annals of Statistics}, 27\penalty0 (3):\penalty0
  1041--1060, 1999.

\bibitem[Tropp(2015)]{MAL-048}
Tropp, J.~A.
\newblock An introduction to matrix concentration inequalities.
\newblock \emph{Foundations and Trends® in Machine Learning}, 8:\penalty0
  1--230, 2015.

\bibitem[Wainwright(2019)]{wainwright2019high}
Wainwright, M.~J.
\newblock \emph{High-dimensional statistics: A non-asymptotic viewpoint},
  volume~48.
\newblock Cambridge University Press, 2019.

\bibitem[Yan et~al.(2012)Yan, Yang, and Xu]{yan2012sparse}
Yan, T., Yang, Y., and Xu, J.
\newblock Sparse paired comparisons in the bradley-terry model.
\newblock \emph{Statistica Sinica}, pp.\  1305--1318, 2012.

\bibitem[Zermelo(1929)]{Zermelo1929}
Zermelo, E.
\newblock Die {B}erechnung der {T}urnier-{E}rgebnisse als ein {M}aximumproblem
  der {W}ahrscheinlichkeitsrechnung.
\newblock \emph{Math. Z.}, 29\penalty0 (1):\penalty0 436--460, 1929.

\end{thebibliography}

\onecolumn
\newpage \appendix

\section{Proofs} \label{app:proofs}

In the following proofs, $n_{ij}$ denotes the number of comparisons between $i$ and $j$. Accordingly, the Laplacian of the comparison graph consists of elements
\begin{equation} \label{eq:laplacian_by_x_ij}
\begin{split}
    L_{ij}
    &= \begin{cases}
    -\frac{n_{ij}}{n}, & i \neq j, \\
    -\sum_{k: k \neq i} L_{ik}, & i = j.
    \end{cases}
\end{split}
\end{equation}
Denoting by $p_{ij}$ the winning probability of $i$ against $j$ for $i, j = 1, \dots, d$, we model the number of comparisons where $i$ defeated $j$ by a Binomial random variable $A_{ij} \sim \mathrm{Binom}(n_{ij}, p_{ij})$. Then, by straightforward calculations,
\begin{equation} \label{eq:ll_by_x_ij}
    l(w) := - \frac{1}{n}\sum_{i,j: i>j} \left\{
        A_{ij} \log(1 + e^{w_j - w_i})
        + A_{ji} \log(1 + e^{w_i - w_j})
    \right\},
\end{equation}
\begin{equation} \label{eq:score_by_x_ij}
    [\nabla l(w^*)]_i = \frac{1}{n}\sum_{j:j \neq i} \left(A_{ij} - n_{ij}p_{ij}\right),
\end{equation}
and the Fisher information matrix (the Hessian of the negative log-likelihood) at $w^*$ has $(i,j)$ coordinates of the form 
\begin{equation} \label{eq:FI_by_x_ij}
\begin{split}
    [\mathcal{I}(w^*)]_{ij}
    &= \begin{cases}
    -\frac{n_{ij}p_{ij}p_{ji}}{n}, & i \neq j, \\
    -\sum_{k: k \neq i} [\mathcal{I}(w^*)]_{ik}, & i = j
    \end{cases}
\end{split}
\end{equation}
for $i, j \in [d]$. 

\subsection{Proof for \texorpdfstring{\cref{thm:existence}}{Theorem 3.2}} \label{app:pf_existence}
For each non-empty $I \subset [d]$, we denote by $E_I$ the event that none of the items in $I^c := [d] \setminus I$ has defeated any item in $I$. That is, 
\begin{equation}
    E_I = \{A_{ij} = n_{ij}, \forall i \in I, j \in I^c\}.
\end{equation}
Since $n_{ij}$ is the maximum value that each $A_{ij}$ can take, the event $E_I$ can be equivalently expressed as
\begin{equation}
    E_I = \left\{\sum_{i \in I, j \in I^c} A_{ij} = \sum_{i \in I, j \in I^c} n_{ij}\right\}.
\end{equation}
Hence, 
\begin{equation}
\begin{split}
    \Pr[E_I] 
    &= \Pr\left[\sum_{i \in I, j \in I^c} A_{ij} = \sum_{i \in I, j \in I^C} n_{ij}\right] \\
    &= \Pr\left[\sum_{i \in I, j \in I^c} (A_{ij} - n_{ij}p_{ij}) = \sum_{i \in I, j \in I^C} n_{ij}(1-p_{ij})\right] \\
    &\leq \Pr\left[\sum_{i \in I, j \in I^c} (A_{ij} - n_{ij}p_{ij}) \geq \sum_{i \in I, j \in I^C} n_{ij}p_{ij}(1-p_{ij})\right].
\end{split}
\end{equation}
Because the $A_{ij}$'s are independent Binomial random variables, we use  Bernstein's inequality to control the probability of the event $E_I$:
\begin{equation}
\begin{split}
    \Pr[E_I] 
    &\leq \Pr\left[\sum_{i \in I, j \in I^c} (A_{ij} - n_{ij}p_{ij})
    \geq \sum_{i \in I, j \in I^c} n_{ij}p_{ij}(1-p_{ij})\right] \\
    &\leq \exp\left(-\frac{3}{2}\sum_{i \in I, j \in I^c} n_{ij}p_{ij}(1-p_{ij})\right).
\end{split}
\end{equation}
The above  bound  can be expressed in terms of the Fisher information matrix $\mathcal{I}(w^*)$. To see this, by the formula \cref{eq:FI_by_x_ij} of the Fisher information matrix, we have that, for each $i \in I$ and $j \in I^c$,  $n_{ij}p_{ij}(1-p_{ij})  = - n \mathcal{I}_{ij}(w^*)$. Next, letting $e_I$ be the indicator vector of the subset $I$ (i.e., the $i$-th element of $e_I$ is $1$ for $i \in I$ and $0$ elsewhere) and noting that $e_{I^c }= {\bf 1} - e_{I}$, we obtain that
\[
\sum_{i \in I, j \in I^c} n_{ij}p_{ij}(1-p_{ij}) = - e_{I} \mathcal{I}(w^*) e_{I^c} =  - e_{I} \mathcal{I}(w^*) ({\bf 1} - e_{I}) = e_{I} \mathcal{I}(w^*) e_{I}, 
\]
where the last identity follows the fact that ${\bf 1}$ is in the null space of $\mathcal{I}(w^*)$. This is because the comparison graph is connected and therefore $\mathcal{I}(w^*)$  has one zero eigenvalue with $\mathrm{span}(\mathbf{1})$ as the associated eigenspace. Thus, we  conclude that
\begin{equation}
\Pr[E_I] \leq \exp\left(-\frac{3n}{2} e_I^\top \mathcal{I}(w^*) e_{I}\right).
\end{equation}
Next, using again the fact that the null space of $\mathcal{I}(w^*)$ is spanned by ${\bf 1}$, 
\begin{equation}
    e_I^\top \mathcal{I}(w^*) e_I 
    \geq \lambda_2(\mathcal{I}(w^*)) \norm{\mathrm{proj}_{\mathrm{span}(\mathbf{1})^\perp}(e_I)}_2^2
    = \lambda_2(\mathcal{I}(w^*)) \frac{\abs{I}(d-\abs{I})}{d},
\end{equation}
where $\mathrm{proj}_{\mathrm{span}(\mathbf{1})^\perp}(\cdot)$ is the  projection mapping onto the linear subspace of $\mathbb{R}^d$ orthogonal to $\mathrm{span}(\mathbf{1})$. 

In the last step of the proof, we follow the arguments from Lemma 1 in \cite{simons1999asymptotics}. By the union bound, 
\begin{equation}
\begin{split}
    \Pr[\text{\cref{cond:existence} fails}] 
    &\leq \sum_{I \in 2^{[d]}} \exp\left(-\frac{3n}{2d}\lambda_2(\mathcal{I}(w^*))\abs{I}(d-\abs{I})\right) \\
    & \leq 2 \sum_{i=1}^{\lceil d/2 \rceil} \binom{d}{i} \exp\left(-\frac{3n}{2d}\lambda_2(\mathcal{I}(w^*))i(d-i)\right) \\
    & \leq 2 \sum_{i=1}^{\lceil d/2 \rceil} \binom{d}{i} \exp\left(-\frac{3n}{2d}\lambda_2(\mathcal{I}(w^*))i\frac{d}{2}\right) \\
    & \leq 2 \left[\left(1+\exp\left(-\frac{3n}{4} \lambda_2(\mathcal{I}(w^*))\right)\right)^d - 1\right].
\end{split}
\end{equation}
The desired result in \Cref{eq:existence} is a straightforward result after plugging-in $\lambda_2(\mathcal{I}^*) \geq 2\frac{\log d}{n}$.

\subsection{Proof for \texorpdfstring{\cref{thm:consistency}}{Theorem 4.1}} \label{app:pf_consistency}

Let $l^*(w)$ be the expected log-likelihood function at $w \in \reals^d$ under the groundtruth score parameter $w^*$, i.e.
\begin{equation}
    l^*(w) := - \frac{1}{n}\sum_{i,j: i>j} \left\{
        n_{ij}p_{ij}  \log(1 + e^{w_j - w_i})
        + n_{ij}p_{ji} \log(1 + e^{w_i - w_j}).
    \right\}
\end{equation}
Since $\log(1+e^{w_j - w_i}) = (w_j - w_i) + \log(1+e^{w_i - w_j})$, we have that
\begin{equation}
    l^*(w) := - \frac{1}{n}\sum_{i,j: i>j} \left\{
        n_{ij}p_{ij}  (w_j - w_i)
        + n_{ij} \log(1 + e^{w_i - w_j})
    \right\},
\end{equation}
and
\begin{equation}
    l(w) := - \frac{1}{n}\sum_{i,j: i>j} \left\{
        A_{ij} (w_j - w_i)
        + n_{ij} \log(1 + e^{w_i - w_j})
    \right\}.
\end{equation}
Thus,
\begin{equation}
    \Delta l(w) := l(w) - l^*(w) = - \frac{1}{n} \sum_{i,j: i > j} \left(A_{ij} - n_{ij}p_{ij}\right) (w_j - w_i) = w^\top \nabla l(w^*).
\end{equation}

Now, let $W$ denote the set of parameter values for which the log-likelihood function is higher than at the true score parameter, i.e., $W := \{w \in \mathbb{R}^d : {\bf 1}^\top w = 0, \;\; l(w) \geq l(w^*)\}$. Since the MLE is one of the elements of $W$, if $W$ is $\ell_2$-bounded by some $r > 0$, then the MLE also exists and is $\ell_2$-bounded by $r > 0$. Throughout the rest of this proof, we let $\hat{w}$ be an arbitrary element in $W$ (possibly the MLE), and $\Delta := \hat{w} - w^*$.
Due to the higher log-likelihood condition of $W$,
\begin{equation}
\begin{split}
    l^*(w^*) - l^*(\hat{w}) &\leq (l^*(w^*) - l^*(\hat{w})) + (l(\hat{w}) - l(w^*)) \\
    &\leq \Delta l(\hat{w}) - \Delta l(w^*) = \Delta^\top \nabla l(w^*).
\end{split}
\end{equation}
On the other hand,
\begin{equation}
\begin{split}
    &l^*(w^*) - l^*(\hat{w}) \\
    &= \frac{1}{n}\sum_{i,j: i>j} \left\{
        n_{ij}p_{ij}  \log\left(\frac{1 + e^{\hat{w}_j - \hat{w}_i}}{1 + e^{w^*_j - w^*_i}}\right) 
        + n_{ij}p_{ji} \log\left(\frac{1 + e^{\hat{w}_i - \hat{w}_j}}{1 + e^{w^*_i - w^*_j}}\right)
    \right\}  \\
    &= \frac{1}{n} \sum_{i,j: i > j} n_{ij} p_{ij} p_{ji} f(w^*_i - w^*_j, \Delta_i - \Delta_j) \\
    &= \sum_{i,j: i > j} -[\mathcal{I}(w^*)]_{ij} f(w^*_i - w^*_j, \Delta_i - \Delta_j)
\end{split}
\end{equation}
where
\begin{equation}
    f(x,y) = (1+e^{-x}) \log\left(\frac{1+e^{x+y}}{1+e^x}\right)
    + (1+e^{x})\log\left(\frac{1+e^{-x-y}}{1+e^{-x}}\right).
\end{equation}

We observe that the curvature of $f(w^*_i - w^*_j, \Delta_i - \Delta_j)$ with respect to $\Delta$ converges to $0$ as $\Delta \rightarrow \infty$, and this property prevents $l^*(w^*) - l^*(\hat{w})$ from being strongly convex. We bypass this problem by deploying the proxy function $h_d$ as defined in \cref{sec:consistency} and show the strong convexity of $l^*(w^*) - l*(\hat{w})$ and $f(w^*_i - w^*_j, \Delta_i - \Delta_j)$ with respect to $h_d(\Delta)$. First, we show in the next result that  $f(w^*_i - w^*_j, \Delta_i - \Delta_j)$ is lower bounded by a quadratic function of $h_d(\Delta_i)$; see \cref{app:pf_f_geq_h^2} for the proof.
\begin{lemma} \label{lem:f_geq_h^2}
    $f(w^*_i - w^*_j, \Delta_i - \Delta_j) \geq c_a (h(\Delta_i) - h(\Delta_j))^2$
    for some universal constant $c_a > 0$.
\end{lemma}
Thus,
\begin{equation}
\begin{split}
    l^*(w^*) - l^*(\hat{w})
    &= \sum_{i,j: i > j} -[\mathcal{I}(w^*)]_{ij} f(w^*_i - w^*_j, \Delta_i - \Delta_j) \\
    &\geq c_a h_d(\Delta)^\top \mathcal{I}(w^*) h_d(\Delta),
\end{split}
\end{equation}
which implies strong convexity of $l^*(w^*) - l^*(\hat{w})$ with respect to $h_d(\Delta)$.
Using this, we arrive at the key basic inequality
\begin{equation}\label{eq:basic_ineq}
    \begin{split}
        \Delta^\top \nabla l(w^*) \geq 
        c_a h_d(\Delta)^\top \mathcal{I}(w^*) h_d(\Delta).
    \end{split}
\end{equation}

We prove \cref{thm:consistency} starting from the basic inequality and further relying on the following, easily verifiable but  special property of $h$: 
\begin{equation} \label{eq:h^2=id-2*h}
    \sgn(\Delta_i) h(\Delta_i)^2 = \Delta_i - 2 h(\Delta_i), ~~\text{for}~~ i = 1, \dots, d.
\end{equation} 
From this, we arrive at the following, nontrivial result; see \cref{app:pf_lemmas} for the proof.
\begin{lemma} \label{lem:rms(h)_geq_2*am^2(h)}
    Provided that $\sum_i \Delta_i = 0$,
    \begin{equation}
        {\sum}_i (h(\Delta_i))^2 \geq \frac{2}{d} \left({\sum}_i h(\Delta_i)\right)^2.
    \end{equation}
\end{lemma}
With the help of the lemma, we obtain that 
\begin{equation} \label{eq:h^2_leq_2*(h-11h)^2}
    \norm{h_d(\Delta)}_2^2 
    \leq 2\norm{h_d(\Delta)}_2^2 - \frac{2}{d} \left({\sum}_i h(\Delta_i)\right)^2 
    = 2 h_d(\Delta)^\top \left(I - \frac{1}{d} \mathbf{1}\mathbf{1}^\top \right) h_d(\Delta).
\end{equation}

Plugging \cref{eq:h^2=id-2*h} into the basic inequality \cref{eq:basic_ineq}, we get that
\begin{equation}
    h_d(\Delta)^\top \mathcal{I}(w^*) h_d(\Delta) 
    \leq  c_{a} \left\{ 2h_d(\Delta)^\top \nabla l(w^*) +
    \sum_{i=1}^d \sgn(\Delta_i) h(\Delta_i)^2 \nabla_i l(w^*) \right\}.
\end{equation}

Hence, 
\begin{equation} \label{eq:quad_ineq_alt}
\begin{split}
    h_d(\Delta)^\top \mathcal{I}(w^*) h_d(\Delta) 
    & \leq c_{a} \left\{2h_d(\Delta)^\top \nabla l(w^*) +
    \norm{\nabla l(w^*)}_\infty \norm{h_d(\Delta)}_2^2 \right\} \\
    & \leq 2 c_{a} \left\{h_d(\Delta)^\top \nabla l_d(w^*) +
    \norm{\nabla l(w^*)}_\infty h_d(\Delta)^\top \left(I - \frac{1}{d} \mathbf{1}\mathbf{1}^\top \right) h_d(\Delta) \right\},
\end{split}
\end{equation}
where the last step is an application of \cref{eq:h^2_leq_2*(h-11h)^2}. Moving the quadratic terms to the left-hand side, 
\begin{equation}
    h_d(\Delta)^\top \left(\mathcal{I}(w^*) - 2c_{a} \norm{\nabla l(w^*)}_\infty \left(I - \frac{1}{d}\mathbf{1}\mathbf{1}^\top\right)\right) h_d(\Delta) \leq 
    2 c_{a} h_d(\Delta)^\top \nabla l(w^*).
\end{equation} 
Now, suppose that it holds that
\begin{equation} \label{eq:I_succeq_dl}
    \mathcal{I}(w^*) \succeq 4 c_a \norm{\nabla l(w^*)}_\infty \left(I - \frac{1}{d}\mathbf{1}\mathbf{1}^\top\right)
\end{equation} 
with high probability;  we will prove this fact below. Then, taking this as a given, we have that
\begin{equation}
     \frac{1}{2} h_d(\Delta)^\top \mathcal{I}(w^*) h_d(\Delta) \leq 
    2 c_{a} h_d(\Delta)^\top \nabla l,
\end{equation} 
and, by Lemma 9 of \citet{shah2016estimation}, 
\begin{equation}
     \frac{1}{2} \norm{h_d(\Delta)}_{\mathcal{I}(w^*)}^2 \leq 
    2 c_{a} \norm{h_d(\Delta)}_{\mathcal{I}(w^*)} \norm{\nabla l(w^*)}_{\mathcal{I}(w^*)^{+}}.
\end{equation} 
As a result,
\begin{equation}
    \norm{h_d(\Delta)}_{\mathcal{I}(w^*)} \leq 4 c_{a} \norm{\nabla l(w^*)}_{\mathcal{I}(w^*)^{+}} \leq 4 c_{a} \sqrt{\kappa} \norm{\nabla l(w^*)}_{L^{+}},
\end{equation}
where the last step is straightforward from the definition of $\kappa$ in \cref{thm:consistency}. Using the Hanson-Wright inequality \citep{hanson1971bound} as in the proof of Theorem 1 of \cite{shah2016estimation}, we obtain that
\[
\norm{\nabla l(w^*)}_{L^{+}}^2 \leq \frac{td}{n},
\]
with probability at least $1 - e^{-t}$. Putting everything together, we thus finally have that
\begin{equation}
    \norm{h_d(\Delta)}_{\mathcal{I}(w^*)}^2 \leq c_1 \kappa \frac{td}{n}
\end{equation}
with probability at least $1 - e^{-t} - p$ where $p$ is the probability that \cref{eq:I_succeq_dl} does not hold. 

To complete the proof, it remains to prove that \cref{eq:I_succeq_dl} holds with high probability. To that effect, we have that, for eah $i \in [d]$,
\begin{equation}
    [\nabla l(w^*)]_i = \frac{1}{n}\sum_{j:j \neq i} \left(x_{ij} - \frac{x_{ij} + x_{ji}}{1+e^{w^*_j - w^*_i}}\right)
\end{equation}
is a sub-Gaussian random variable with $\sigma^2 = v_{\max}/4$. By a standard maximal inequality for sub-Gaussian random variables \citep{wainwright2019high}, 
\begin{equation}\label{eq:max_ineq_1}
    \Pr\left[\norm{\nabla l(w^*)}_\infty \geq \frac{\sqrt{v_{\max}}}{2n} (\sqrt{2 \log(2d)} + t) \right] \leq e^{-t^2/2}
\end{equation}
for any $t > 0$. Setting $t = \sqrt{2 \log(2d)}$, we have that
\begin{equation}\label{eq:max_ineq_2}
    \Pr\left[\norm{\nabla l(w^*)}_\infty \geq \frac{1}{n}\sqrt{2 v_{\max} \log(2d)}\right] \leq \frac{1}{2d}.
\end{equation}
Hence, provided that $\lambda_2(\mathcal{I}(w^*)) \geq \frac{c_0}{n} \sqrt{v_\text{max} \log d}$ for some large enough $c_0 > 0$, $\lambda_2(\mathcal{I}(w^*)) \geq \norm{\nabla l(w^*)}_\infty$ with probability at least $1 - \frac{1}{2d}$, for sufficiently large $d$. This follows by a straightforward derivation of \cref{eq:I_succeq_dl}. 

\subsection{Proofs of Lemmas} \label{app:pf_lemmas}

In this section we provide the proofs of technical Lemmas needed in the proof of \cref{thm:consistency}.

\subsubsection{Proof of \texorpdfstring{\cref{lem:lambda2L}}{Lemma 4.2}}
\label{app:pf_lambda2L}

In this proof, we use that this graph has a similar algebraic connectivity to the Cayley graph on $(\ints_d, +)$ with difference set $\{-W, \dots,-1\} \cup \{1, \dots, W\}$. This is formalized in the following argument, which leverages results form the spectral graph theory for the Cayley graph (e.g., \citealp{brouwer2011spectra})

Let $\tilde{L}$ be the normalized graph Laplacian of the Cayley graph on $(\ints_d,+)$ with difference set  $\{-W, \dots -1\} \cup \{1, \dots, W\}$. $\tilde{L}$ has elements
\begin{equation}\label{eq:tilde_L}
    \tilde{L}_{ij} = \begin{cases}
    -\sum_{k: k \neq i} \tilde{L}_{ik}, & i = j, \\
    -1/\tilde{n}, & 0 < \abs{i - j} \leq W, \\
    -1/\tilde{n}, & 0 < d - \abs{i - j} \leq W, \\
    0, & \text{elsewhere,}
    \end{cases}
\end{equation}
where $\tilde{n} = dW$. It follows from Proposition 1.7.1 in \citet{brouwer2011spectra} that $\lambda_2(\tilde{L}) \geq \lambda_2(L)$. 

On the other hand, let $\bar{L}$ be the normalized graph Laplacian of the Cayley graph on $(\ints_{2d},+)$ with the same difference set. $\bar{L}$ can be obtained from \cref{eq:tilde_L} by changing $d$ to $2d$. Suppose that we fold this graph so that vertices $i$ and $2d-i$ are equated. The obtained folded graph has normalized graph Laplacian $L'$ such that
\begin{equation}\label{eq:L'}
    L'_{ij} = \begin{cases}
    -\sum_{k: k \neq i} L'_{ik}, & i = j, \\
    -3/n', & i \neq j \textand i+j \leq W+2, \\
    -3/n', & i \neq j \textand i+j \geq d-W, \\
    0, & \abs{i - j} > W, \\
    -1/n', & \text{elsewhere,}
    \end{cases}
\end{equation}
where $n' = 2dW$. If $(x_1, \dots, x_d)$ is an eigenvector of $L'$ with eigenvalue $\lambda$, then $(x_1, \dots, x_d, x_d, \dots, x_1)$ is an eigenvector of $\bar{L}$ with the same eigenvalue. Hence, $\lambda_2(\bar{L})$ is always smaller than $\lambda_2(L')$. Now, we notice that $L'_{ij} \leq 3 L_{ij}$ for every $(i,j) \in [d]^2$ and hence Proposition 1.7.1 in \citet{brouwer2011spectra} again implies $\lambda_2(L') \leq 3 \lambda_2(L)$. Putting everything together,
\begin{equation}
    \frac{1}{3} \lambda_2(\bar{L}) \leq \frac{1}{3} \lambda_2(L') \leq \lambda_2(L) \leq \lambda_2(\tilde{L}).
\end{equation}

According to the spectral properties of Cayley graphs listed in \citet{brouwer2011spectra}, 
\begin{equation}
\begin{split}
    \lambda_2(\tilde{L}) &\approx \frac{2W}{n} - \frac{2}{n}  \sum_{i=1}^W \cos\left(\frac{2\pi i}{d}\right) \\
    &= \frac{2}{n} \sum_{i=1}^W \left(1-\cos\left(\frac{2\pi i}{d}\right)\right) \\
    &= \Theta\left(\frac{d}{n\pi} \int_0^{2\pi W/d} (1-\cos x) dx\right)\\
    &= \Theta\left(\frac{W^3}{nd^2}\right).
\end{split}
\end{equation}
Similarly, $\lambda_2(\bar{L}) = \Theta\left(\frac{W^2}{d^3}\right)$, and $\lambda_2(L)$ has the same rate in $d$.


\subsubsection{Proof of \texorpdfstring{\cref{lem:f_geq_h^2}}{Lemma A.1}}
\label{app:pf_f_geq_h^2}
Let
\begin{equation}
    g(y) := \log\left(\frac{1+e^y}{2}\right) + \log\left(\frac{1+e^{-y}}{2}\right).
\end{equation}

We will derive the desired inequality from the following two properties of $f$, $g$, and $h$:
\begin{itemize}
    \item $f(x,y) \geq g(y)$ for any $x , y \in \reals$, and
    \item $g(y_1-y_2) \geq c_a \left(h(y_1) - h(y_2)\right)^2, \forall y_1, y_2 \in \reals$ for some universal constant $c_a$.
\end{itemize}

First, for any $x, y \in \reals$,
\begin{equation}
    \frac{\partial}{\partial y} f(x,y) = (1+e^{-x}) \frac{e^{x+y}}{1+e^{x+y}} - (1+e^x) \frac{e^{-x-y}}{1+e^{-x-y}}
\end{equation}
and
\begin{equation}
    \frac{\partial}{\partial y} g(y) = \frac{e^{y}}{1+e^{y}} - \frac{e^{-y}}{1+e^{-y}}.
\end{equation}
Hence,
\begin{equation}
\begin{split}
    \frac{\partial}{\partial y} \left(f(x,y) - g(y)\right) &= \frac{e^x+e^y}{(1+e^{x+y})(1+e^y)(1+e^{-y})} (e^y - e^{-y}) \\
    &\begin{cases}
    > 0, & y > 0, \\
    = 0, & y = 0, \\
    < 0, & y < 0.
    \end{cases}
\end{split}
\end{equation}
With the fact that $\frac{\partial}{\partial y} \left(f(x,y) - g(y))\right) \big|_{y=0} = 0$, we derive the first inequality between $f$ and $g$.

Now, we move on to the second inequality between $g$ and $h$. We first note that
\begin{equation}
    \sqrt{x+h} - \sqrt{x} \geq \sqrt{y+h} - \sqrt{y}
\end{equation}
where $0 < x < y$ and $h > 0$, which can be easily seen through the following inequality:
\begin{equation}
    \sqrt{x+h} - \sqrt{x} = \frac{h}{\sqrt{x+h}+\sqrt{x}} \geq \frac{h}{\sqrt{y+h}+\sqrt{y}} = \sqrt{y+h} - \sqrt{y}.
\end{equation}
Hence, 
\begin{equation}
    (h(y_1) - h(y_2))^2 = (\sqrt{\abs{y_1}+1} - \sqrt{\abs{y_2}+1})^2 \leq (\sqrt{\abs{\abs{y_1}-\abs{y_2}}+1} - 1)^2 \leq h(y_1 - y_2)^2.
\end{equation}
That is, it suffices to show that $g(x) \geq c_a h(x)^2$, $\forall x$ for some universal constant $c_a$; then,
\begin{equation}
    c_a (h(y_1) - h(y_2))^2 \leq c_a h(y_1 - y_2)^2 \leq g(y_1 - y_2).
\end{equation}
We note that both $g(\cdot)$ and $h(\cdot)^2$ increase linearly around $0$ and so do $\sqrt{g(\cdot)}$ and $h(\cdot)$ sufficiently away from $0$. In other words, $\frac{h(x)^2}{g(x)}$ converges as either $x \rightarrow 0$ or $\rightarrow \infty$. Being continuous, $\frac{h(x)^2}{g(x)}$ is upper bounded, and we obtain the desired inequality directly.


\section{Supplementary Simulation Results on Different Graph Topologies} \label{app:supp_simulation}

\begin{figure}[t!]
  \centering 
  \setlength{\labelsep}{0mm}
  \sidesubfloat[]{
    \includegraphics[height = 0.24\textwidth]{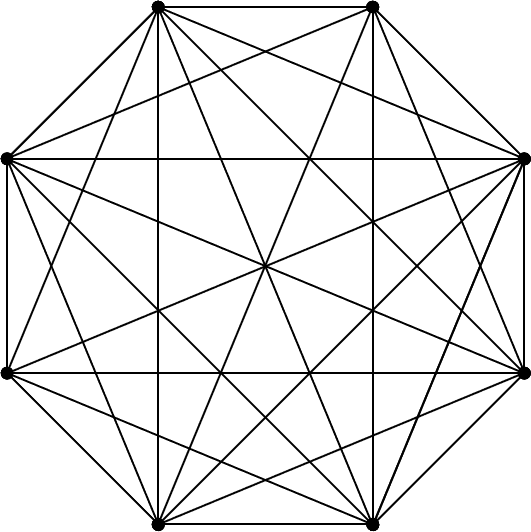}
    \label{fig:graph_complete}
  }
  \setlength{\labelsep}{0mm}
  \sidesubfloat[]{
    \includegraphics[height = 0.24\textwidth]{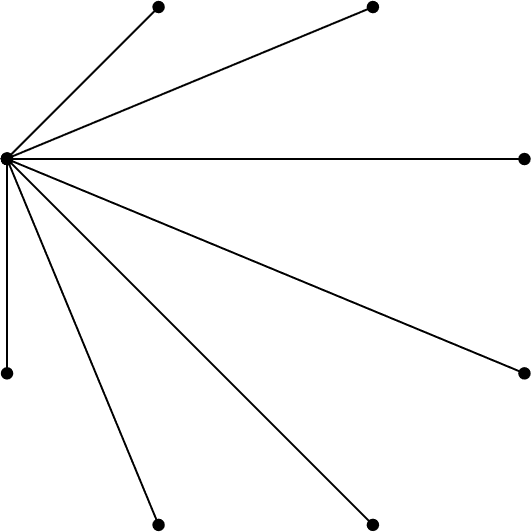}
    \label{fig:graph_star}
  }
  \sidesubfloat[]{
    \includegraphics[height = 0.24\textwidth]{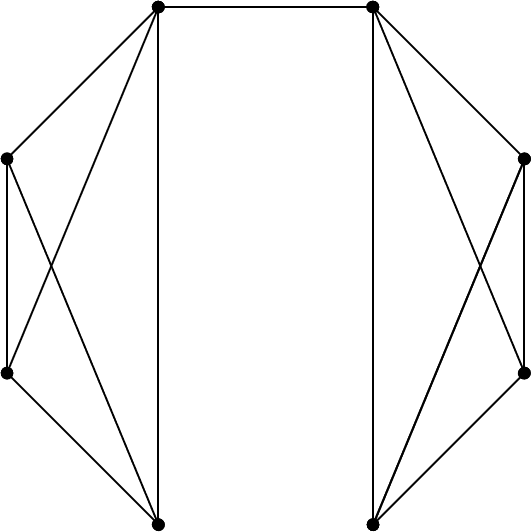}
    \label{fig:graph_barbell}
  }

  \caption{\sl {\bf Illustrations of the three simulated comparison graphs with $d = 8$.} {\bf (a)} complete graph, {\bf (b)} star graph, and {\bf (c)} barbell graph. 
  }
  \label{fig:other_graphs}
\end{figure}

\begin{figure}[t!]
  \centering 
  \setlength{\labelsep}{-4mm}
  \sidesubfloat[]{
    \includegraphics[height = 0.25\textwidth]{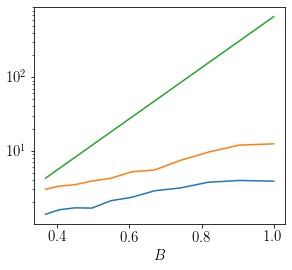}
    \label{fig:l2_vs_B_complete}
  }
  \setlength{\labelsep}{-4mm}
  \sidesubfloat[]{
    \includegraphics[height = 0.25\textwidth]{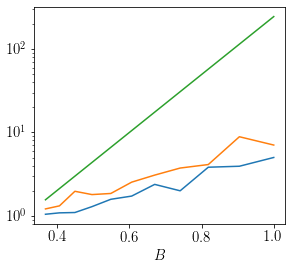}
    \label{fig:l2_vs_B_star}
  }
  \sidesubfloat[]{
    \includegraphics[height = 0.25\textwidth]{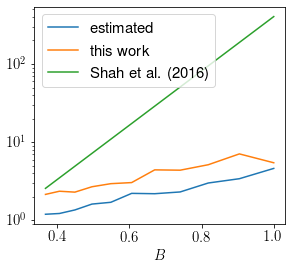}
    \label{fig:l2_vs_B_barbell}
  }

  \caption{\sl {\bf $\ell_2$ error versus the dynamic range parameter $B$.} Average empirical $\ell_2$ estimation error evaluated via simulations for {\bf (a)} complete graphs, {\bf (b)} star graphs, and {\bf (c)} barbell graphs. The value of $B$ ranges from $e^{-1}$ to $1$ for all three graphs, where $d$ was set to $100$ for the complete graph and to $20$ for the other graphs. The figures also depict the rate of increase of the theoretical upper bounds based on  \cref{thm:consistency} and Theorem 2 of \citet{shah2016estimation}. 
  }
  \label{fig:l2_vs_B_other_graphs}
\end{figure}

\begin{figure}[t!]
  \centering 
  \setlength{\labelsep}{-1mm}
  \sidesubfloat[]{
    \includegraphics[height = 0.25\textwidth]{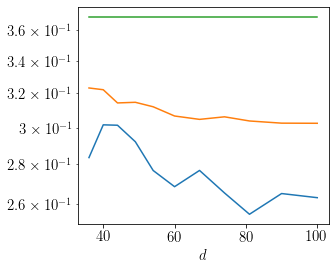}
    \label{fig:l2_vs_d_complete}
  }
  \setlength{\labelsep}{-1mm}
  \sidesubfloat[]{
    \includegraphics[height = 0.25\textwidth]{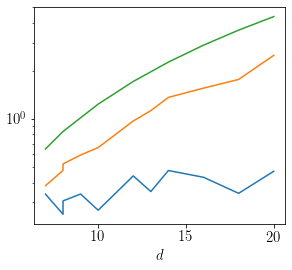}
    \label{fig:l2_vs_d_star}
  }
  \sidesubfloat[]{
    \includegraphics[height = 0.25\textwidth]{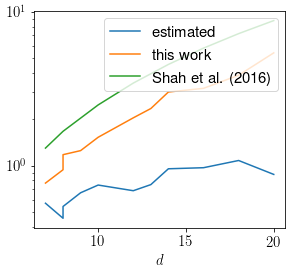}
    \label{fig:l2_vs_d_barbell}
  }

  \caption{\sl {\bf $\ell_2$ error versus  the number of items $d$.} Average empirical $\ell_2$ estimation error evaluated via simulations for {\bf (a)} complete graphs, {\bf (b)} star graphs, and {\bf (c)} barbell graphs. For complete graphs, $d$ ranges from $\floor{100e^{-1}}$ to $100$, where the range is from $\floor{20e^{-1}}$ to $20$ for the others. $B$ was fixed at $1$ for every case. The figures also depict the rate of increase of the theoretical upper bounds based on  \cref{thm:consistency} and Theorem 2 of \citet{shah2016estimation}. 
  }
  \label{fig:l2_vs_d_other_graphs}
\end{figure}

We conducted supplementary simulation experiments for the three comparison graph topologies studied by simulation in \citet{shah2016estimation}, which are listed below.
Unfortunately, there does not exist a parameter modulating the connectivity for the three graph topologies, such as the comparison width $W$ in the banded comparison graphs. We instead evaluated via simulations the $\ell_2$ error bound of the maximum likelihood BTL parameter estimator as a function of the dynamic range parameter $B$ and the number of items $d$ and compared to the rate of increase associated to the theoretical bounds from \cref{thm:consistency} and Theorem 2, \citet{shah2016estimation}.
Similarly to the simulation settings in \citet{shah2016estimation}, the true score parameters $w^*$ were obtained as a normalized sample from the standard $d$-dimensional Gaussian distribution: for each $i =1 ,\ldots, d$, we sampled $w'_i \distiid \distNorm(0,1)$ and then set $w^*_i = B \cdot {w'_i}/{\norm{w'}_\infty}$ to ensure a dynamic range value of $B$.
The three kinds of comparison graphs and simulation settings for this study are as follow: 
\begin{enumerate}
    \item {\it Complete graph:} every pair of $d$ items are compared $T=5$ times. $B$ ranges from $e^{-1}$ to $1$, and $d$  from $\floor{100 e^{-1}}$ to $100$.
    \item {\it Star graph:} one {\it hub} item is compared $T=200$ times against every other item, and there is no comparison among the $d-1$ non-hub items. $B$ ranges from $e^{-1}$ to $1$, and $d$ from $\floor{20 e^{-1}}$ to $20$.
    \item {\it Barbell graph:} for even $d$, every pair of items $1, \dots, d/2$ are compared $T=100$ times, and so are items $d/2+1, \dots, d$. Items $d/2$ and $d/2+1$ are also compared $T=100$ times. $B$ ranges from $e^{-1}$ to $1$, and $d$  from $\floor{20 e^{-1}}$ to $20$.
\end{enumerate}
See \cref{fig:other_graphs} for the illustrations of the simulated comparison graphs when $d = 8$.
For each of the three graph topologies and  a grid of values for $B$ and $d$, we simulated $100$ datasets, each time with a different draw for $w^*$. The maximum likelihood BTL parameter estimate for each dataset was found as described in \cref{sec:simulation}.
Since a different value of $w^*$ was used for each dataset, a direct comparison between the average estimated $\ell_2$ errors and  the theoretical error bounds can be misleading. Instead, we evaluate  the 95th percentile of the empirical $\ell_2$ errors over the simulated datasets for each combination of $(B,d)$.
\cref{fig:l2_vs_B_other_graphs} and \cref{fig:l2_vs_d_other_graphs} show the 95th percentile curves of the  empirical $\ell_2$ errors for the three graph typologies shown in \cref{fig:other_graphs} as a function of the dynamic range $B$ and of the number of items $d$, respectively. On each plot we also depict the theoretical bounds from \cref{thm:consistency} and Theorem 2 of \citet{shah2016estimation} arbitrarily shifted for better readability. Thus,  the curves are to be compared based on their rates of growth of and not their values. The plot demonstrate that, for a fixed $d$ and varying $B$, the error rates predicted by \cref{thm:consistency} track closely the estimated $\ell_2$ errors and are more accurate than the ones implied by Theorem 2 of \citet{shah2016estimation} in all the graph topologies under consideration. On the other hand, the two error bounds appear to be similar and possibly loose as a function of $d$ for a fixed value of $B$.
\cref{fig:l2_vs_B_other_graphs} shows  

.

\end{document}